\documentclass{article}
\usepackage{hyperref}
\hypersetup{
  colorlinks   = true, 
  urlcolor     = black, 
  linkcolor    = blue, 
  citecolor   = green 
}
\usepackage{geometry}
\usepackage[english]{babel}
\usepackage{amsmath}
\usepackage{graphicx}
\usepackage{ marvosym }
\usepackage[utf8]{inputenc}
\usepackage{mathtools}
\usepackage{geometry}
\usepackage{amsfonts}
\usepackage{amssymb}
\usepackage{amsthm}
\usepackage{thmtools}
\usepackage{t1enc}
\usepackage[titles]{tocloft}
\usepackage{makeidx}
\usepackage{wasysym}
\usepackage{stmaryrd}
\usepackage{calc}  
\usepackage{enumitem} 
\usepackage[refpage]{nomencl}

\usepackage{algpseudocode}
\usepackage{tikz}
\usetikzlibrary{arrows}
\usepackage{float}

\DeclarePairedDelimiterX\Set[2]{\lbrace}{\rbrace}%
 { #1 \,\delimsize:\, #2 }

\theoremstyle{definition}

\newenvironment{biz}{\par\noindent{\itshape Proof:}\ }{\rule{1.5ex}{1.5ex}}
\newenvironment{sbiz}{\par\noindent{\itshape Proof:}\ }{\newmoon}

\theoremstyle{plain}
\newtheorem{thm}{Theorem}

\newtheorem{claim}[thm]{Claim}
\newtheorem{obs}[thm]{Observation}
\newtheorem{prop}[thm]{Proposition}
\newtheorem*{prop*}{Proposition}
\newtheorem*{seged*}{Sublemma}

\newtheorem{lem}[thm]{Lemma}

\newtheorem*{cond*}{Condition}

\newtheorem*{lem*}{Lemma}
\theoremstyle{definition}

\newtheorem*{defn*}{Definition}

\newtheorem{fel*}[thm]{Exercise}

\newtheorem*{megf*}{Observation}
\theoremstyle{remark}
\newtheorem{rem}[thm]{Remark}
\newtheorem*{rem*}{Remark}

\title{T-joins in infinite graphs as edge-disjoint system of paths matching the vertices in $ T $}
\author{Attila Joó\thanks{MTA-ELTE Egerváry Research Group, Department of Operations Research, Eötvös University, Budapest, Hungary. 
E-mail: {\tt 
joapaat@cs.elte.hu}.}}
\date{2016}

\begin{document}

\maketitle

\begin{abstract}
 We characterize the class of infinite connected graphs $ G $ for which there  exists a $ T $-join for any choice of an infinite $ T 
 \subseteq V(G) $. We also show that the following well-known fact remains true in the infinite case. If $ G $ is connected and does not 
 contain a $ T  $-join, then it will if we either remove an arbitrary vertex from $ T $ or add any new vertex to $ T $.
\end{abstract}

\section{Introduction}
The graphs in this paper may have multiple edges although all of our results follow easily from their restrictions to simple graphs. 
Loops are irrelevant, hence we throw them away automatically if during some construction arise some. The $ 2 $-edge-connected components 
of a graph are its maximal $ 2 $-edge-connected subgraphs (a graph consists of a single vertex is considered $ 2 $-edge-connected).  A $ 
\boldsymbol{T} $\textbf{-join} in a graph $ G $, where $ 
T\subseteq V(G) $, is a system $ \mathcal{P} $ of edge-disjoint paths in $ G $ such that the endvertices of the 
paths in $ \mathcal{P} $ create a partition of $ T $ into two-element sets. In other words, we match by edge-disjoint paths the 
vertices in $ T $. In the finite case the existence of an $ F\subseteq E(G) $ for which $ d_F(v) $ is odd if and only if $ v\in T $ is 
equivalent 
with the existence of a $ T $-join. Indeed, the united edge sets of the paths in $ \mathcal{P} $ forms such an $ F $, and such an $ F $ 
can 
be decomposed into a $ T $-join and some cycles by the greedy method. Sometimes  $ F $ itself is called a $ T $-join. The two possible 
definitions are no more 
closely related in the infinite case. Take  for example a one-way infinite path where $ T $  contains only its endvertex. 
Then there is no $ T $-join according to the path-system based definition (which we will use during this paper) but the whole edge set 
forms a $ T $-join with respect to the second definition.
 
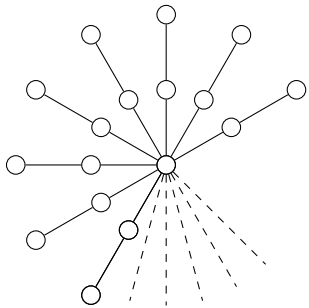
\begin{figure}[H]
 \centering
 
 \begin{tikzpicture}
 
 \node[circle,inner sep=0pt,draw, minimum size=7] (v0) at (0,0) {};
 \node[circle,inner sep=0pt,draw, minimum size=7] (v1) at ({30}:2cm) {};
  \node[circle,inner sep=0pt,draw, minimum size=7] (v2) at ({60}:2cm) {};
  \node[circle,inner sep=0pt,draw, minimum size=7] (v3) at ({90}:2cm) {};
  \node[circle,inner sep=0pt,draw, minimum size=7] (v4) at ({120}:2cm) {};
  \node[circle,inner sep=0pt,draw, minimum size=7] (v5) at ({150}:2cm) {};
  \node[circle,inner sep=0pt,draw, minimum size=7] (v6) at ({180}:2cm) {};
  \node[circle,inner sep=0pt,draw, minimum size=7] (v7) at ({210}:2cm) {};
  \node[circle,inner sep=0pt,draw, minimum size=7] (v8) at ({240}:2cm) {};
  \node[circle,inner sep=0pt,draw, minimum size=7] (v9) at ({240}:2cm) {};

  \node[circle,inner sep=0pt,draw, minimum size=7] (u0) at (0,0) {};
   \node[circle,inner sep=0pt,draw, minimum size=7] (u1) at ({30}:1cm) {};
    \node[circle,inner sep=0pt,draw, minimum size=7] (u2) at ({60}:1cm) {};
    \node[circle,inner sep=0pt,draw, minimum size=7] (u3) at ({90}:1cm) {};
    \node[circle,inner sep=0pt,draw, minimum size=7] (u4) at ({120}:1cm) {};
    \node[circle,inner sep=0pt,draw, minimum size=7] (u5) at ({150}:1cm) {};
    \node[circle,inner sep=0pt,draw, minimum size=7] (u6) at ({180}:1cm) {};
    \node[circle,inner sep=0pt,draw, minimum size=7] (u7) at ({210}:1cm) {};
    \node[circle,inner sep=0pt,draw, minimum size=7] (u8) at ({240}:1cm) {};
    \node[circle,inner sep=0pt,draw, minimum size=7] (u9) at ({240}:1cm) {};

  \node[circle,inner sep=0pt, minimum size=7] (v_10) at ({255}:2cm) {};
  \node[circle,inner sep=0pt, minimum size=7] (v_11) at ({270}:2cm) {};
  \node[circle,inner sep=0pt, minimum size=7] (v_12) at ({285}:2cm) {};
  \node[circle,inner sep=0pt, minimum size=7] (v_13) at ({300}:2cm) {};
  \node[circle,inner sep=0pt, minimum size=7] (v_14) at ({315}:2cm) {};

  \draw  (v0) edge (u1);
  \draw  (v0) edge (u2);
  \draw  (v0) edge (u3);
  \draw  (v0) edge (u4);
  \draw  (v0) edge (u5);
  \draw  (v0) edge (u6);
  \draw  (v0) edge (u7);
  \draw  (v0) edge (u8);
  \draw  (v0) edge (u9);
  
  \draw  (u1) edge (v1);
    \draw  (u2) edge (v2);
    \draw  (u3) edge (v3);
    \draw  (u4) edge (v4);
    \draw  (u5) edge (v5);
    \draw  (u6) edge (v6);
    \draw  (u7) edge (v7);
    \draw  (u8) edge (v8);
    \draw  (u9) edge (v9);
  
  \draw  (v0) edge[dashed] (v_10);
  \draw  (v0) edge[dashed] (v_11);
  \draw  (v0) edge[dashed] (v_12);
  \draw  (v0) edge[dashed] (v_13);
  \draw  (v0) edge[dashed] (v_14);
 
 \end{tikzpicture}
 \caption{ A subdivided infinite star. It has no $ T $-join if $ T $ is the whole vertex set.}
 \label{nincs T-kötés} 
 
 \end{figure}

$ T $-join is a common tool in combinatorial optimization problems such as the well-known Chinese postman 
problem. For a 
detailed survey one can see \cite{as1994survey}. For finite connected graphs the necessary and sufficient condition for the existence of 
a $ T $-join is quiet simple, $ \left|T\right| $ must be even. Indeed, the necessity of the condition is trivial. For the sufficiency 
 let $ \left|T\right|=2k $ and we apply induction on $ k $. The case $ k=0 $ is clear. If $ \left|T\right|=2k+2 $, then remove two 
 vertices, $ u $ and $ v $ say, of $ T $ to obtain $ T' $. By induction we have a $ T' $-join. Take the symmetric difference of the edge 
 set of a $ T' $-join and the edges of an arbitrary path between $ u $ and $ v $. By the greedy method we can partition the resulting 
 edge set into a $ T $-join and some cycles.  If $ \left|T\right| $ is even but $ G $ is 
 infinite, then the same proof works. In this paper we investigate questions related to the existence of $ T $-joins where 
$ T $ is infinite. For infinite $ T $ one 
can not guarantee in general the existence of a $ T $-join in a connected graph. 
Consider for example an infinite star and subdivide all of its edges by a new vertex (see Figure \ref{nincs T-kötés}) and  let $ T $ be 
consist of the whole vertex set. 

 One of Our main results is that essentially Figure \ref{nincs T-kötés} is the only counterexample.
 
 \begin{thm}
 The connected  graph $ G $ does not contain  a $ T $-join for some infinite $ T\subseteq V(G) $ if and only if  one can obtain the 
 subdivided infinite star  
 (see Figure \ref{nincs T-kötés}) from $ G $ by contracting edges and deleting the resulting loops.
 \end{thm}
The following reformulation of the theorem will be more convenient:

\begin{thm}\label{minden végtelenre Tköt}
A connected infinite graph $ G $ contains a $ T $-join for every infinite $ T\subseteq V(G) $ if and only if there is no $ U \subseteq 
V(G) $ for 
which every connected component of $ G-U $ connects to $ U $ in $ G $ by a single edge and infinitely many of them are nontrivial (i.e. 
not consist of a single vertex).
\end{thm}

The ``if'' direction is straightforward since if such an $ U $ exists, then one can choose one vertex from $ U $ and two from each 
nontrivial connected component of $ G-U $  to obtain a $ T $ 
for which there is no $ T $-join in $ G $.

Our other result describes the effect of finite modifications of $ T $ on the existence of a $ T $-join. For finite $ T $  if $ G $ does 
not contain a $ T $-join, then it 
contains a $ T' $-join whenever $\left|(T\setminus T')\cup (T'\setminus T)\right|=: \left|T \triangle T'\right| $ is 
odd. (It is obvious, since in this case the existence of a $ T $-join depends just on the parity of $ \left|T\right| $.) Surprisingly 
this property  
remains true for infinite $ T $ as well. We have the following result about this.
\begin{thm}\label{Tköt lokális viselkedés}
The class  $ \{ (G,T):\ G \text{ is a connnected graph and }T\subseteq V(G) \} $ can be partitioned into three nonempty subclasses 
defined by the following three properties.  

\begin{enumerate}
 \item[(A)]  $ G $ contains a $ T' $-join whenever $ \left|T'\triangle T\right|<\infty $,
 \item[(B)] $ G $ contains  $ T' $-join if $ \left|T'\triangle T\right| $ is finite and even, but it does not when it is odd,
\item[(C)]   $ G $ contains  $ T' $-join if $ \left|T'\triangle T\right| $ is finite and odd, but it does not when it is even.
\end{enumerate} 
\end{thm}

For countable $ T $ our proofs based on purely combinatorial arguments. To handle uncountable $ T $ we apply the so called elementary 
submodel technique. We do not assume previous knowledge about this method. We build it up shortly and recommend 
\cite{soukup2011elementary} for a more detailed introduction. 

To make the descriptions of the $ T $-join-constructing processes more reader-friendly we introduce the following 
single player game terminology. 
There is an abstract set of tokens and every token is on some vertex of a graph $ G $. (At the beginning typically we have exactly one  
token on each  
element of a prescribed  vertex set $ T $ and none on the other vertices.) If two tokens are on the same vertex, then we may remove them 
(we say that we \textbf{match} them to each other). 
If we have a token $ t $ on the vertex $ u $ and $ uv $ is an edge of $ G $, then we may move $ t $ from $ u $ to $ v $, but then we have 
to delete $ uv $ from $ G $. A \textbf{gameplay} is a transfinite sequence of the steps above in which we move every token just finitely 
many 
times. Limit steps are defined by the earlier steps in a natural way. Indeed, we just delete all the edges that have been deleted 
earlier, 
remove the 
tokens that have been removed before, and put all the remaining tokens to their stabilized positions. We call a gameplay winning if we 
remove 
all the tokens eventually. To make  talking about the relevant part of the graph easier we also allow as a feasible step  deleting 
vertices  
without tokens on them and deleting edges. Clearly  there is a $ T $-join in $ G $ if and only if  there is a  winning gameplay for the 
game 
on $ G $ 
where the 
initial token distribution is 
 defined by $ T $.

\section{The 2-edge-connected case}
A subgraph of $ G $ is called $ \boldsymbol{t} $\textbf{-infinite} if it contains infinitely many tokens. We define the notions $ 
\boldsymbol{t} 
$\textbf{-finite}, $ \boldsymbol{t} $\textbf{-empty}, $ \boldsymbol{t} $\textbf{-odd} and $ \boldsymbol{t} $\textbf{-even} similarly.
\begin{lem}\label{Tkötés 2élöfre}
If $ G $ is 2-edge-connected and $ t $-infinite, then there is a winning gameplay.
\end{lem}
\begin{claim}\label{tokenek elkerülő matchelése}
Assume that $ G $ is 2-edge-connected and contains even number of tokens, but at least four including $ s\neq t $.  Then there 
is a winning gameplay in which $ s $ and $ t $ are not matched with each other and $ t $ is not moved.
\end{claim}
 \begin{sbiz}
 We may assume that $ G $ is finite otherwise we may take a finite 2-edge-connected subgraph that contains all the tokens. Take two 
 edge-disjoint paths between the vertices that contain  $ s $ and $ t $ and let $ H $ be the Eulerian subgraph of $ G $ consists of these 
 paths. We can build 
 up $ G $ from $ H $ by adding ears (as in the ear decomposition). We apply induction on the number of ears. If there is no ear i.e. $ 
 G=H $, then we take an Eulerian cycle $ O $ in $ G $.  Fix an Eulerian orientation of $ O $. Either this orientation or the reverse of 
 it induces a desired gameplay.
 
  Otherwise 
 let $ Q $ be the last ear. If the number of tokens on the new vertices given by $ Q $ is odd (even), then match inside $ Q $ all but one 
 (two) of these tokens  and 
 move the exception(s)  to the part of the $ G $ before the addition of $ Q $. Delete the remaining part of $ Q $. We are done by 
 applying the induction hypothesis.   
 \end{sbiz}\\

By contracting the $ 2 $-edge-connected components of a connected graph, we obtain a tree. If $ R $ is a 2-edge-connected component of a 
connected graph $ G $, then we denote by $ \boldsymbol{\mathsf{tree}(G; R)} $ the tree of 
the 2-edge-connected components of $ G $ rooted at $ R $. We usually pick such a root $ R $ arbitrarily without mentioning it 
explicitly.  
We write $ \mathsf{tree}(G) $ if it is not rooted. We do not distinguish 
strictly the subtrees of $ \mathsf{tree}(G; 
R) $ and the corresponding subgraphs of $ G $.

\begin{claim}\label{rekonstruál alak(T-kötés)}
Let  $ G $ be a  connected graph with infinitely many tokens on it. Assume that $ G $ has only finitely many 2-edge-connected 
components.  Then there is 
a finite gameplay after which all the components of the resulting $ G' $ are 2-edge-connected and contain infinitely many 
tokens. 
\end{claim}
\begin{sbiz}
By the pigeon hole principle there is a $ t $-infinite 2-edge-connected component $ R $ of $ G $. We use induction on the number $ k $ 
of the 
2-edge-connected components. For $ k=1 $ we do nothing. If $ k>1 $ we take a leaf $ C $ of $ \mathsf{tree}(G;R) $. If $ C $ is $ t 
$-infinite we remove the 
unique outgoing edge of $ V(C) $ in $ G $ and we use induction to the  arising component other than $ C $. If $ C $ is $ t 
$-even, 
then we match the tokens on $ C $ with each other inside $ C $ and we delete the remaining part of $ C $ and its unique outgoing edge and 
use induction. In 
the $ t $-odd case we match all but one tokens of $ C $ inside $ C $ and move one to the parent component. It is doable by adding another 
``phantom-token'' $ t $ on the vertex of $ C $ incident with
the cut edge to the parent component and applying Claim \ref{tokenek elkerülő matchelése} with this $ t $ and an arbitrary $ s $. We 
delete the remaining part 
of $ C $ again and use induction.
\end{sbiz}\\

 Now we turn to the proof of Lemma \ref{Tkötés 2élöfre}. If $ t $ is a token and $ H $ is a 
  subgraph of $ G $, then we use the abbreviation 
  $ \boldsymbol{t\in H} $ to express the fact that $ t $ is on some vertex of $ H $. Assume first that $ T $ is 
   countable. For technical reasons we assume just the following weaker condition instead of $ 2 $-edge-connectedness.
 
 \begin{equation}\label{Tj componensek 2ec és t-inf}
 \text{All the connected components are 2-edge-connected and t-infinite.}
 \end{equation}
 
   Let  $ t_0 $ be an arbitrary token.  It is enough to show, that there is a finite gameplay such that we remove 
 $ t_0 $ and  the resulting system still satisfies the condition \ref{Tj componensek 2ec és t-inf}.  Pick two edge-disjoint paths $ 
 P_1,P_2 $
between $ t_0 $ and any other token $ t^{*} $ that lies in the same component as $ t_0 $. For $ i\in \{ 1,2 \} $ let $ K_i $ be the 
connected component of $ G-E(P_i) $ that contains $ t^{*} $ and $ t_0 $. We claim that either $ K_1 $ or $ K_2 $ is $ t $-infinite. 
Suppose that $ 
K_1 $ is not. Then there is a $ t $-infinite component $ K_{inf} $ of $ G-E(P_1) $ which does not contain $ t^{*} $ and $ t_0 $. Note 
that $ P_2 $ 
necessarily lies in $ K_1 $. Hence $ K_{inf} $ is a subgraph of $ G-E(P_2) $ and $ P_1 $ ensures that $ K_{inf} $ and $ K_1 $ belongs to 
the same connected component of $ G-E(P_2) $. Thus $ K_{inf} $ ensures that this component is $ t $-infinite. We will need the following 
basic observation.

\begin{obs}\label{Tj nem lesz sok 2ec}
If each component of a graph $ G $ has finitely many
$ 2 $-edge-connected components, then so has $ G-f $ for every $ f\in E(G) $.
\end{obs} 

By symmetry we may assume that $ K_1 $ is $ t $-infinite. Move $ t_0 $ along the edges of $ P_1 $ one by one. If the following edge $ e $ 
is a bridge, and moving $ t_0 $ along $ e $ would create a $ t $-odd component $ K_{odd} $, then $ t^{*}\notin K_{odd} $ because the 
component which contains  $ t^{*} $  is $ t $-infinite. In this case we delete $ e $ without moving $ t_0 $ and obtain 
a  $ t $-infinite  component and a $ t $-even component $ K_{even} $ that contains $ t_0 $. We match the tokens on $ K_{even} $ and erase 
the remaining part of it and the first phase of the process is 
done. If this case does not occur, then we move $ t_0 $ to $ t^{*} $ 
along $ P_1 $ and remove both.   

We need to fix the condition \ref{Tj componensek 2ec és t-inf}. Each component of the resulting graph is either $ t $-even or $ t 
$-infinite. We match the tokens on $ t $-even components and erase the 
remaining part of them. Each $ t $-infinite component has finitely many $ 2 $-edge-connected components by Observation 
\ref{Tj nem lesz sok 2ec} thus we are done by applying Claim \ref{rekonstruál alak(T-kötés)}.

Consider now the general case where $ T $ can be arbitrary large. Examples show that our approach for countable $ T $ may fail for 
uncountable $ T $ at limit steps. Add a new vertex $ 
z $ and draw all the edges $ zv\ (v\in 
T) $  to 
obtain $ H $. Finding a $ T $-join in $ G $  is obviously equivalent with covering in $ H $ all the edges incident with $ z $ by 
edge-disjoint cycles. To reduce the problem to the countable case it is enough to prove the following claim.

\begin{claim}\label{Tkötés elemi részmodell}
There is a partition of $ E(H) $ into countable sets $ E_i\ (i\in I) $ in such a way, that for all $ i\in I $  the graphs $ 
H_i:=(V(H),E_i) 
$ 
have the following property. The connected components of $ H_i-z $ are 2-edge-connected and connect to $ z $ in $ H_i $ by either zero or 
infinitely 
many 
edges.
\end{claim}

Our proof of Claim \ref{Tkötés elemi részmodell} is a basic application of the  elementary submodel technique.  One 
can find a detailed survey about this method with many combinatorial applications in \cite{soukup2011elementary}. We give here 
just the fundamental definitions and cite the results that we need.
   Let $ 
\Sigma=\{ \varphi_1,\dots, \varphi_n \} $ be a finite set of formulas in the language of set theory where the free variables of $ 
\varphi_i $ are $ 
x_{i,1},\dots,x_{i,n_i} $. A set $ M $ is  a $ \boldsymbol{\Sigma} $\textbf{-elementary submodel} if the formulas in $ \Sigma $ are 
absolute between $ (M,\in) $ and the universe i.e.  
\[ \left[(M,\in) \models \varphi_i(a_1,\dots,a_{n_i}) \right] \Longleftrightarrow  \varphi_i(a_1,\dots,a_{n_i}) \] 
\noindent holds whenever $ 1\leq i \leq n $ and $ a_1,\dots,  a_{n_i} \in M $.  By using Lévy's Reflection Principle and the 
Downward 
Löwenheim Skolem Theorem (see in \cite{kunen2011set} or any other set theory or logic textbook) one can derive the following fact.

\begin{claim}\label{Tköt van elemirészmodell}
For  any finite se $ \Sigma $  of formulas,  set $ x $ and infinite cardinal $ \kappa $ there exists a $ \Sigma $-elementary submodel 
$ M\ni x $ 
with $\kappa= \left|M\right|\subseteq M $.
\end{claim} 

Now we use some  methods developed by  L. Soukup in \cite{soukup2011elementary}. Call a class $  \mathfrak{C}  $ of graphs
\textbf{well-reflecting} if for all large enough finite set $ \Sigma $ of formulas, infinite  cardinal $ \kappa $, set $ x $ and $ G\in 
\mathfrak{C} $ there is a  $ \Sigma $-elementary submodel $ M $ with $x,G\in M $ for which  $ \kappa=\left|M\right|\subseteq M $ and 
$ (V(G), E(G)\cap M),\  (V(G),E(G)\setminus M)\in \mathfrak{C} $. (``For all large enough finite $ \Sigma $'' means here that there is 
some finite $ \Sigma_0 $ such that for all finite $ \Sigma \supseteq \Sigma_0$.)

\begin{thm}[L. Soukup, Theorem 5.4 of \cite{soukup2011elementary}]
Assume that the graph-class $ \mathfrak{C} $ is well-reflecting and $ G\in \mathfrak{C} $. Then there is a partition of $ E(G) $ into 
countable 
sets $ E_i\ (i\in I) $ in such a way, that for all $ i\in I $ we have $ (V(G),E_i)\in \mathfrak{C} $.
\end{thm}

\begin{rem}
L. Soukup used originally a stricter notion of well-reflectingness but his proof still works with our weaker notion as well. 
\end{rem}

 We apply the Theorem above to prove Claim \ref{Tkötés elemi részmodell}. Let $ \mathfrak{C} $ be the class of graphs $ G $ for 
 which $ z\in V(G) $, the connected components of $ G-z $ are 2-edge-connected and connect to $ z $ in $ G $ either by 
 infinitely many edges or send no edge to $ z $ at all. We need to show that $ \mathfrak{C} $ is well-reflecting. Assume that $ \Sigma $ 
 is a finite set of formulas that contains all the formulas of length at most $ 10^{10} $ with variables at most $ x_1,\dots x_{10^{10}} 
 $. (From the proof one can get an exact list of formulas need to be in $ \Sigma $. The usual terminology says just to fix a large enough 
 $ 
 \Sigma $. We decided that a more explicit definition is beneficial for readers who first met whit this method.) 
 
 Let  $ \kappa,x$ and $ 
 G\in \mathfrak{C} $ be given. By Claim \ref{Tköt van elemirészmodell} we can find a $ \Sigma $-elementary submodel $ M \ni 
 x,G,z $ with $ \kappa= \left|M\right| \subseteq M $. We know that  $ (G-z)\in M $ by using the 
absoluteness of the formula  ``$ x_1 $ graph obtained by the deletion of vertex $ x_2 $ of graph $ x_3 $''$ \in \Sigma $ (see Claim 2.7 
and 2.8 of \cite{soukup2011elementary} for 
basic facts about absoluteness). The proof of  $ 
(V(G),E(G)\cap M)\in \mathfrak{C} $ is easy. We just need the absoluteness of formulas such that ``the local edge-connectivity between 
the 
vertices $ x_1 $ and $ x_2 $ in the graph $ x_3 $ is $ x_4 $'' $ \in \Sigma $.
 The hard part is to show $ 
(V(G),E(G)\setminus M)\in \mathfrak{C} $. We use the following proposition.
\begin{prop}[Lemma 5.3 of \cite{soukup2011elementary}]\label{egeszben szeparal}
 If $ M $ is a $ \Sigma $-elementary submodel (for some large enough finite $ \Sigma $) for 
which $ G\in M $, $ \left|M\right|\subseteq M $ and
 $ x\neq y\in V(G) $ are in the same connected component of $ (V(G),E(G)\setminus M) $ and $ F\subseteq E(G) \setminus M $ separates them 
 where $ \left|F\right| \leq \left|M\right|$; then $ F $ separates $ x $ and $ y $ in  $ G $ as well.
\end{prop}

 If the local edge-connectivity between vertices $ x\neq y $ would be one in the graph  $ (V(G-z),E(G-z)\setminus M) $, then  we can 
 separate 
 them by the deletion of a single edge $ e $. But then by applying Proposition \ref{egeszben szeparal} with $ 
F=\{ e \} $ we may conclude that the same separation is possible in $ G-z $, which contradicts  the assumption $ G\in \mathfrak{C} $.
 
Suppose, to the contrary, that $ z $ sends finitely many, but at least one, edges, say $ e_1,\dots,e_k $,  to a 
connected component of $ (V(G-z),E(G-z)\setminus M) $ in $ (V(G),E(G)\setminus M) $. Let $ p $ be the endvertex of $ e_1 $ other than $ z 
$. Then $ F:=\{ e_i 
\}_{i=1}^{k} $ separates $ p $ and $ z $ in $  
(V(G),E(G)\setminus M)  $ and  $ \left|F\right|<\infty $  holds.  Hence $ F $ separates them in $ G $ as well which is a contradiction. 
Now the proof of Lemma \ref{Tkötés 2élöfre} is complete.

\section{The simplification process}

Lemma \ref{Tkötés 2élöfre} makes it possible to decide the existence of a $ T $-join by just investigating the structure of the 
2-edge-connected components and the quantity of tokens on them.   Let $ G $ be a connected graph, $ T \subseteq V(G) $ and let $ R $ be a 
2-edge-connected 
component of it.  We define a graph $ H=H(G,R,T) $ with a token-distribution on it. To obtain $ H $ we apply the following gameplay that 
we call 
\textbf{simplification process}. We denote  by $ \boldsymbol{\mathsf{subt}(C;G,R)} $ the  subtree of the descendants of  $ C $ rooted at 
$ C $  in $ \mathsf{tree}(G;R) $.   
Delete all those 2-edge-connected 
componets  $ 
C $ for which $ \mathsf{subt}(C;G,R) $ does not contain 
any token.   Then consider 
the $ t $-finite leafs  
of the reminder of $ 
\mathsf{tree}(G;R) $. (We do not consider the root as a leaf.) Match the tokens on any $ t $-even 
leaf  $ C $  inside $ C $ and for $ t $-odd leafs $ C $ move one token to the parent and match the others inside. In both cases 
delete the 
remaining part of $ C $. Iterate the steps above as long as possible and denote by $ H $ the 
resulting graph. Clearly either $ H=R $ 
with some tokens on it or if $ \mathsf{tree}(H) $ is nontrivial, then $ \mathsf{subt}(C;H, R) $ must be $ t 
$-infinite 
for all 2-edge-connected components $ C $ of $ H $.

\begin{claim}\label{Tköt t-végtelen megállás nyerő}
There is a winning gameplay for the original system if and only if $ H $ is $ t $-even or $ t $-infinite.
\end{claim} 

\begin{sbiz}
We show the ``if'' part here and the ``only if part'' later in Claim \ref{Tköt nem vesztünk megoldást}. In the $ t $-even case it is 
obvious since $ H $ is connected. 
Assume that $ H $ is $ t $-infinite. We may suppose that $ H $ is not 2-edge-connected because otherwise we are done by applying Lemma 
\ref{Tkötés 2élöfre}. 
Then $ 
\mathsf{subt}(C;H,R) $ is $ t 
 $-infinite for all 2-edge-connected components $ C $ of $ H $.  For each  2 -edge-connected 
  components of $ H $ we  
 fix a path $ P_C $ in the tree $ \mathsf{subt}(C;H,R) $ that starts at $ C $ and either terminates at some leaf 
 of $ \mathsf{tree}(H;R) $ or it is 
 one-way infinite and meets some 2-edge-connected  not $ t $-empty component  of $ H $ other than $ C $.  
 
 After these 
 preparations we do the following. If the root $ R $ is $ t $-even or $ t $-infinite, then we match 
all the tokens of it inside $ R $ (use Lemma \ref{Tkötés 2élöfre} in the $ t $-infinite case) and delete the remaining part of $ R $. If 
it is $ t $-odd we move one of its tokens, say it will turn to be $ t^{*} $, to some child  of $ R $ determined by the path $ P_R $ and 
we define $ P_{t^{*}}:=P_R $. We match the 
other tokens inside $ R $ and then delete the remaining part of $ 
R 
$. At the next step we 
deal with the $ \mathsf{subt}(C;H,R) $ trees where $ C $ is a child of $ R $. In the cases where $ C $ is  $ t $-infinite or $ t $-even 
we 
do 
the same as 
earlier. Assume that $ C $ is $ t $-odd. If there is no token on $ C $ that comes from $ R $, then we do the same as earlier. Suppose 
that there is, say $ t_0 $. If there is a token on $ C $ other than $ t_0 $, then we match here $ t_0 $ and send forward some other 
token $ t_1 $ in the direction defined by $ P_C $ (apply Claim \ref{tokenek elkerülő matchelése} and a phantom-token) and we let $ 
P_{t_1}:=P_C $. If $ t_0 $ is the only token of $ C_0 $, then we move $ 
t_0 $ in the direction $ P_{t_0} $. We iterate the process recursively. The only not entirely trivial thing that we need to justify is 
that we do not move a token infinitely many times. If we moved some token $ t $ at the previous step, then we match it at the current 
step unless it is the only token at the corresponding 2-edge-connected component. On the other hand, when we move $ t $ for the first 
time we define the path $ P_t $. The later movements of $ t $ are leaded by $ P_t $ which ensure that eventually $ t $ will meet 
some other token.     
\end{sbiz}

\section{Proof of the theorems}

Now we are able to prove the nontrivial direction of Theorem \ref{minden végtelenre Tköt}. \begin{biz}
Let $ G $ be an infinite connected graph 
 such that there is no $ U \subseteq V(G) $ for 
 which the connected components of $ G-U $ connect to $ U $ in $ G $ by a single edge and infinitely many of them are nontrivial. Let $ 
 T \subseteq V(G) $ be 
infinite. Consider 
the vertices $ v $ of degree one that are in $ T $ i.e. there is a token on them. Move these tokens to the only possible direction and 
then delete all the vertices of degree one or zero. We denote the resulting graph by $ G' $ and we fix a  2-edge-connected 
component $ R $ of it. If $ G' $ has a $ 
t $-infinite 2-edge-connected 
component, then  it cannot vanish during the simplification process, thus the resulting $ H $ will be $ t $-infinite and we are done by 
 Claim \ref{Tköt t-végtelen megállás nyerő}.  

Assume there is no such  
component. 
 Degree of $ C_0:=R $ in $ \mathsf{tree}(G') $ must be finite 
otherwise  $ U:=V(C_0) $ would violate the assumption about $ G $. Since $ C_0 $ is $ t $-finite by the pigeonhole principle there is a 
child $ C_1 $ of $ C_0 $ such 
that $ \mathsf{subt}(C_1; G', R) $ contains infinitely many tokens. By recursion we obtain a one-way infinite path of $ \mathsf{tree}(G') 
$ with 
vertices $ C_n\  (n\in 
\mathbb{N})$  such that for every $ n $ the tree $ \mathsf{subt}(C_n; G', R) $ contains infinitely many tokens. The set $ 
U:=\bigcup_{n=0}^{\infty}V(C_n) $ may not have infinitely many outgoing edges in $ G' $ otherwise $ U $  violates the condition about $ G 
$. Thus for large enough $ n $ the tree $ \mathsf{subt}(C_n; G', R) $ is just a terminal segment of the one-way infinite path we 
constructed. 
This implies 
that infinitely many of the $ 
C_n $'s 
contain at least one token. Since such a path cannot vanish during the simplification process, it terminates with a $ t $-infinite $ H $ 
again. \end{biz}\\ 

\noindent We turn to the proof of Theorem \ref{Tköt lokális viselkedés}.
 \begin{biz}
 Let a connected $ G $ and a 2-edge-connected component $ R $ of it be 
fixed. Let $ T 
\subseteq V(G) $. We will 
show that case (A) of Theorem \ref{Tköt lokális viselkedés} occurs if and only if the simplification process terminates with infinitely 
many tokens  
and (B)/(C) 
 occurs if and only if it terminates with an even/odd number of tokens respectively.\\ 
 
 Assume 
first 
that  the result $ H $ of the simplification process initialized by $ T $ ($ T $-process from now on)  is  $ t $-infinite.  Let $ 
T'\subseteq V(G) $ such that $ \left|T' \triangle T\right|<\infty $. Call $ T' 
$-process 
the simplification process with the initial tokens given by $ T' $ and  denote by $ H' $ the result of it. If $ G $ has a 
2-edge-connected, $ t 
$-infinite component $ C $ with respect to $ T $, then $ C $ is $ t $-infinite with respect to $ T' $ as well. Observe that such a $ C $ 
remains untouched during the simplification process. Thus $ H' $ is $ t $-infinite  and 
therefore there is a $ T' $-join in $ G $ by Claim \ref{Tköt t-végtelen megállás nyerő}.

  We may suppose that there is no 2-edge-connected, $ t $-infinite component in $ G $. If such a component $ C $ arises during the $ T 
  $-process, then $ C $ receives a token from infinitely many children  of it. Since $ \left|T \triangle T'\right|<\infty $ we have $ 
  \left|V(D)\cap T\right|=\left|V(D)\cap T'\right| $ for all but finitely many 2-edge-connected component $ D $. Hence the 
  token-structure of 
  $ \mathsf{subt}(D; G, R) $ is the same for all but finitely many child  $ D $ of $ C $ at the  case of $ T' $. Thus $ C $ will receive 
  infinitely many tokens during the $ T' $-process as well. 
  
  Finally 
  we  suppose that such a component does not arise, i.e. $ H $ has no $ t $-infinite component. Then $ 
  \mathsf{tree}(H) $ must be an infinite tree since $ H $ is $ t $-infinite. Furthermore, $ \mathsf{subt}(C; H, R) $ must contain at 
  least 
  one 
  token 
  for all 2-edge-connected component $ C $ of $ H $ otherwise we may erase $ \mathsf{subt}(C;H,R) $ to continue the simplification 
  process. Fix a one-way infinite path $ P $ in $ \mathsf{tree}(H) $ with vertices $ C_n (n\in \mathbb{N}) $, where $ C_0=R $ and 
  infinitely many $ C_n $ contain at least one token. For a large enough $ n_0 $ the token-distribution of $ \mathsf{subt}(C_{n_0};G,R) 
  $ 
  is the same at the $ T $ and at the $ T' $ cases. Hence the $ T $-process and $ T' $-process runs identical on the subgraph of $ G $ 
  corresponding to $ \mathsf{subt}(C_{n_0};G,R) $.  
  Thus   $  \mathsf{subt}(C_{n_0};H, R)=\mathsf{subt}(C_{n_0}; H', R) $ holds and the tokens on them are the same. But then the terminal 
  segment of path  $  P $ in $ \mathsf{tree}(H') $ shows that $ H' $ is $ t $-infinite as well. 

\begin{claim}\label{Tköt nem vesztünk megoldást}
There is no $ T $-join in $ G $ if the simplification process terminates with an odd number of tokens.
\end{claim}

\begin{sbiz}
Remember that no $ t $-infinite 2-edge-connected component  may arise during the simplification process in this case. Assume, to the 
contrary, that 
there is a $ T $-join $ \mathcal{P} $ in $ G $. We play a winning gameplay induced by $ \mathcal{P} 
$ i.e. every step we move some token along the appropriate $ P\in \mathcal{P} $ towards its partner. If for a $ 2 $-edge-connected 
component $ C $  the 
subgraph 
$ \mathsf{subt}(C;G,R) $ does not contain any vertex from $ T $, then clearly no $ P\in \mathcal{P} $ comes here. Hence we may delete 
these 
parts 
of the graph. If $ C $ is a leaf of (the remaining part of) $ 
\mathsf{tree}(G;R)  $ with $ \left|T\cap V(C)\right| $ even, then the corresponding paths are inside $ C $. In the odd case  exactly 
one $ P\in \mathcal{P} $ uses the unique outgoing edge of $ C $ and some other paths match inside $ C $ the other $ T $-vertices of $ C 
$. Thus along $ P $ we may move one token to the parent component and match the others along the other paths in $ C $.   Therefore after 
we do 
these 
steps the quantity of the tokens on the 2-edge-connected components will be the same as after the first step of the simplification 
process.  Similar 
arguments show that it remains true after successor steps of the simplification process as well.  Since the first difference may not 
arise at a limit step for all steps of the simplification process we have a corresponding position of the gameplay induced by $ 
\mathcal{P} 
$ where the token 
quantities on the $ 2 $-edge-connected components are the same. On 
the other hand, it cannot be true for the terminating position since the play induced by $ \mathcal{P} $ is a winning gameplay and hence 
it 
cannot arise a system with an odd number of tokens during the play which is 
a contradiction.
\end{sbiz} 

\begin{claim}\label{Tköt odd even váltás}
If the simplification process  for $ T $ terminates with an even (odd) number of tokens and $ \left|T \triangle T'\right|=1 $, then 
the 
simplification process for  $ T' $ terminates with an odd (even) number of tokens.
\end{claim} 

\begin{sbiz}
  By symmetry we may let  $ 
 T'=T\cup\{ v \} $. If $ v\in V(R) $, then the simplification process for $ T' $ runs in the same way as for $ T $ except 
 that 
 at the end we have the extra token on $ v $ which changes the parity of the remaining tokens as we claimed. If $ v $ is not in the root 
 $ R $, 
 then it is in $ \mathsf{subt}(C;G,R) $ for some child $ C $ of $ R $. This $ C $ is closer  in $ \mathsf{tree}(G) $  to the 
 2-edge-connected component that contains $ v $ than $ R $. On the one hand we know by 
 induction that 
 the parity of the number of tokens on $ C $ will be different when $ C $ will become a leaf  in the case of $ T' $. On 
 the other hand for the other children of $ R $ this parity will be clearly the same, thus the parity of the number of the remaining 
 tokens changed again. 
\end{sbiz}

The remaining part of Theorem \ref{Tköt lokális viselkedés}  follows from  Claim \ref{Tköt nem vesztünk megoldást} and from the 
repeated application of Claim \ref{Tköt odd even váltás}. 
\end{biz}


\begin{thebibliography}{1}

\bibitem{as1994survey}
{\sc Frank, A.}
\newblock A survey on t-joins, t-cuts, and conservative weightings, 1994.

\bibitem{kunen2011set}
{\sc Kunen, K.}
\newblock {\em Set theory}.
\newblock College Publ., 2011.

\bibitem{soukup2011elementary}
{\sc Soukup, L.}
\newblock Elementary submodels in infinite combinatorics.
\newblock {\em Discrete Mathematics 311}, 15 (2011), 1585--1598.

\end{thebibliography}
\end{document}